\documentclass[pdflatex,sn-mathphys-num]{sn-jnl}
\usepackage{graphicx}%
\usepackage{multirow}%
\usepackage{amsmath,amssymb,amsfonts}%
\usepackage{amsthm}%
\usepackage{mathrsfs}%
\usepackage[title]{appendix}%
\usepackage{xcolor}%
\usepackage{textcomp}%
\usepackage{manyfoot}%
\usepackage{booktabs}%
\usepackage{algorithm}%
\usepackage{algorithmicx}%
\usepackage{algpseudocode}%
\usepackage{listings}%
\usepackage{enumitem}%
\usepackage{bm}
\usepackage{enumitem}

\usepackage{booktabs} 
\usepackage{caption}  
\usepackage{array}    
\usepackage{float}  
\theoremstyle{thmstyleone}%
\newtheorem{theorem}{Theorem}[section]
\newtheorem{lemma}[theorem]{Lemma}
\newtheorem{proposition}[theorem]{Proposition}
\newtheorem{corollary}[theorem]{Corollary}
\newtheorem{conjecture}[theorem]{Conjecture}

\theoremstyle{thmstyletwo}%
\newtheorem{example}[theorem]{Example}
\newtheorem{remark}[theorem]{Remark}

\theoremstyle{thmstylethree}%
\newtheorem{definition}[theorem]{Definition}
\newtheorem{open problem}{Open problem}
\begin{document}
	
	\title[Combinatorial $t$-Designs from Finite Abelian Groups and Their Applications to Elliptic Curve Codes]{Combinatorial $t$-Designs from Finite Abelian Groups and Their Applications to Elliptic Curve Codes}

	\author[1]{\fnm{Hengfeng} \sur{Liu}}\email{hengfengliu@163.com}
	
	\author*[2]{\fnm{Chunming} \sur{Tang}}\email{tangchunmingmath@163.com}
	
	\author[3]{\fnm{Cuiling} \sur{Fan}}\email{cuilingfan@163.com}
	
	\author[4]{\fnm{Rong} \sur{Luo}} \email{luorong@swjtu.edu.cn}
	\affil[1,3,4]{\orgdiv{School of Mathematics}, \orgname{Southwest Jiaotong University}, \city{Chengdu}, \postcode{611756},  \country{China}}
	
	\affil[2]{\orgdiv{School of Information Science and Technology}, \orgname{Southwest Jiaotong University}, \city{Chengdu}, \postcode{611756}, \country{China}}


	\abstract{In this paper, we establish the conditions for some finite abelian groups and the family all the $k$-sets in each of them summing up to an element $x$ to form $t$-designs. We fully characterize the sufficient and necessary conditions for the incidence structures to form $1$-designs in finite abelian $p$-groups, generalizing existing results on vector spaces over finite fields. For finite abelian groups of exponent $pq$, we also propose sufficient and necessary conditions for the incidence structures to form a $1$-designs. Furthermore, some interesting observations of the general case when the group is cyclic or non-cyclic are presented and the relations between $(t-1)$-designs and $t$-designs from subset sums are established. As an application, we demonstrate the correspondence between $t$-designs from the minimum-weight codewords in elliptic curve codes and subset-sum designs in their groups of rational points. By such a correspondence, elliptic curve codes supporting designs can be simply derived from subset sums in finite abelian groups that supporting designs.}

	\maketitle

	\section{Introduction}
	Let $v, k, \lambda \in \mathbb{N}$. A $t$-$\left ( v, k, \lambda   \right ) $ \emph{design} is an incidence structure $\mathcal{D}= \left (\mathcal{G} ,\mathcal{B}   \right ) $, where $\mathcal{G}$ is a set with $v$ elements (called points) and $\mathcal{B}$ is a family of distinct subsets  of $\mathcal{G}$ (called blocks), such that any block in $\mathcal{B}$ contains exactly $k$ points of $\mathcal{G}$, and any subset of  $\mathcal{G}$ of size $t$ is contained in exactly $\lambda $ blocks in $\mathcal{B}$. In a $t$-$\left ( v, k, \lambda   \right ) $  design, the parameters are interrelated, for $t=1$, let $b$ denote the number of blocks, then the equation $v\lambda=bk$ is obtained by counting the pairs $(B,p)$, where $B$ is a block and $p$ is a point contained in $B$. Similarly, for $t=2$, a $2$-$\left ( v, k, \lambda   \right ) $ design is also a $1$-$\left ( v, k, r \right ) $ design, with the equation $\lambda (v-1)=r(k-1)$, obtained by counting the triples $(s,t, B)$ where $s$ is a fixed point and $t$ is a distinct point from $s$, and $B$ is a block containing both of them. Generally, a $t$-$\left ( v, k, \lambda_{t}   \right ) $ design is also an $i$-$\left ( v, k, \lambda_{i}  \right ) $ design for $1\le i\le t-1$, where $$\lambda_{i}=\lambda_{t}\binom{v-i}{t-i} \Big/ \binom{k-i}{t-i}.$$ 
	The number of blocks $b_{t}$ satisfy
	$$\binom{n}{t}\lambda _{t}= \binom{k}{t} b_{t}.$$ 
	For more information about $t$-designs, the reader is referred to \cite{BJL,C2008}.\par
	
	In this paper, we mainly focus on the conditions for the family of subsets in an abelian group that sum to a given element to support a $t$-design. Let $G$ be an additively written abelian group, and let $D\subseteq G$ be a finite subset of order $n$. For an element $x\in G$, we denote the number of $k$-subsets of $D$ that sum up to $x$ by $$N(D,k,x)=\#\left\{T\subseteq D{:}\,\#T=k,\sum_{t\in T}t=x\right\}.$$
	The well-known \emph{subset-sum problem} is to determine whether $N(D,k,x)>0$ for some $1\le k \le n$. This NP-complete problem arises from coding theory and cryptography, with applications in the knapsack cryptosystem (for $G=\mathbb{Z}$), the deep hole problem of extended Reed–Solomon codes (for $G=\mathbb{F}_{q}$) \cite{ZCL2016}, and the minimal distance of elliptic curve codes (for $G=E(\mathbb{F}_{q})$) \cite{C2008}. The main challenge of the problem comes from the flexibility in choosing the subset $D$. Despite the fact that the problem is typically challenging, a lot of progress has been made when the subset $D$ has specific algebraic structures, especially when $D=G$. Li and Wan \cite{LW2008} derived a closed form for $N(G,k,x)$, where $G$ is the additive group of a finite field $\mathbb{F}_{q}$ (an elementary abelian $p$-group). In \cite{LW2012} they extended their results to the case where $G$ is a finite abelian group, and a concise proof using character theory was later provided by Kosters in \cite{K2013}. Furthermore, in \cite{K2013} the author shows that $N(G,k,x)$ is nonzero except in certain trivial cases. \par
	When the underlying group $G$ is settled, we denote the family of $k$-subsets that sum up to $x\in G$ by $\mathcal{B} _{k}^{x}$ and represent $N(G,k,x)$ as $b_{k}^{x}$. The question of whether the incidence structure $(G, \mathcal{B} _{k}^{x} )$ supports a $t$-design has been shown to have strong connections to coding theory \cite{D1973}, particularly in the context of Hamming codes \cite{FP2021}, when $G=\mathbb{F}_{2}^{d}$. Moreover, in \cite{FP2021} and \cite{P2023} the authors established sufficient and necessary conditions for $(\mathbb{F}_{p}^{d}, \mathcal{B} _{k}^{x} )$ to be a $1$-design and $2$-design, where $p$ is a prime. While their results are only for an elementary abelian $p$-group (with exponent $p$), that is, $G=\mathbb{Z}_{p}\oplus\cdots\oplus\mathbb{Z}_{p}$, is natural to ask when $(G, \mathcal{B} _{k}^{x} )$ is a $t$-design where $G$ is a general finite abelian group. We generalize the result of \cite{P2023} which concerns finite elementary abelian $p$-groups, by characterizing the conditions for $(G, \mathcal{B} _{k}^{x} )$ to be a $1$-design, where $G$ is any finite abelian $p$-group. Moreover, we propose a conjecture on the absence of $2$-design in non-elementary abelian $p$-groups. As a further step towards general finite abelian groups, we also characterize the conditions for the incidence structure to be a $1$-design, when $G$ is any finite abelian group with exponent $pq$, where $p,q$ are distinct primes. Besides, we also provide additional observations of the $t$-designs from $(G, \mathcal{B} _{k}^{x} )$, where $G$ is a cyclic group or a finite abelian group of non-cyclic type, which are closely related to the elliptic curve codes.\par 
	An $[n,k,d]$ linear code is called an MDS code if it achieves the Singleton bound, i.e., $d=n-k+1$, and it is said to be almost maximum distance separable (AMDS for short) if $d=n-k$. If a code and its dual code are both AMDS, then the code is said to be near MDS (NMDS for short). NMDS codes are of interest because they have many nice applications in combinatorial designs and cryptography \cite{D22020,Heng2023-NMDS,Simos2012,XuG2022,ZhiY2025}. In recent years, many NMDS codes have been constructed \cite{DY2024,Heng2022,Heng2023-NMDS,LiH 2023,XuL2023,XuL2024,Yin2024,ZhiY2025}. Linear codes that supporting designs have attracted significant attention these years \cite{D12020,D22020,D32020,T2020,W2023,Xiang2022,Xiang2023,Yan2024,Yin2024}, as it is an important approach to construct $t$-designs and codes supporting designs may have efficient decoding methods \cite{T2020, Xiang2023, Yan2024}.  \par
	Elliptic curve codes are a special class of algebraic geometry codes constructed from rational points on elliptic curves over finite fields. These codes offer better parameters compared to classical linear codes (such as BCH or Reed-Solomon codes), making them attractive from application standpoint. It is well-known that the rational points on elliptic curves over finite fields forms an finite abelian group, and we characterize the support designs of minimum-weight codewords in elliptic curve codes from $t$-designs held in their rational points groups. Such a correspondence shows the potential of obtaining a large number of NMDS elliptic curve codes supporting $t$-designs from subset sums in the corresponding groups of rational points.
	\par

	Our main contributions in this paper are as follows:
	\begin{itemize}
		\item For any finite abelian $p$-group $G$ with exponent $p^{m}$ $(m\ge1)$, we fully characterize the conditions for $(G, \mathcal{B} _{k}^{x} )$ to be a $1$-design, and propose a conjecture on the case of $2$-design. Notably, our results include the findings on elementary abelian $p$-groups in \cite{P2023} as a special case.
		\item For any finite abelian group $G$ with exponent $pq$ where $p,q$ are distinct primes, we also characterize the conditions for $(G, \mathcal{B} _{k}^{x} )$ to be a $1$-design, as a further step towards general finite abelian groups.
		\item We make some observations on $t$-designs from incidence structure $(G, \mathcal{B} _{k}^{x} )$, where $G$ is a cyclic group or a finite abelian group of non-cyclic type. In particular, we establish the relations between $(t-1)$-designs and $t$-designs from subset sums , which are closely related to the elliptic curve code. 
		\item As an application, we characterize the support designs of minimum-weight codewords in some elliptic curve codes from $t$-designs held in their rational points groups. By such correspondence, we obtain a class of NMDS elliptic curve codes supporting $1$-designs.
	\end{itemize}
	The remaining sections of this paper are arranged as follows. Section \ref{sec pre} presents some notations and notions. In this section, we also recall some important results from the past research, which will be used in subsequent sections. Section \ref{sec p} studies $1$-design of subset sums in any finite abelian $p$-group $G$ of exponent $p^{m}$ and proposes a conjecture on the case of $2$-design. We investigate $1$-design of subset sums in any finite abelian group $G$ of exponent $pq$, where $p$ and $q$ are two distinct primes, in Section \ref{sec pq} . Section \ref{sec cyc} makes some notes on designs held in cyclic groups and non-cyclic abelian groups. Section \ref{elliptic curve codes} discusses the connection between $t$-designs from subset sums in finite abelian groups and some NMDS elliptic curve codes are derived. Section \ref{sec con} concludes the paper with a short summary followed by two open problems.

	\section{Preliminaries} \label{sec pre}
	This section gives closed forms of the coefficients introduced in the following definition and also covers other useful results from the earlier works. Throughout the paper, $\mathbb{F}_q$ denotes the finite field of order $q$, where $q$ is some prime power. Let $G$ be a finite abelian group and let $G^{*}=G\setminus \left \{ 0 \right \} $. 
	\begin{definition}\label{def coee}
		Let $G$ be a finite abelian group of order $n$ and let $k$ be a positive integer satisfying $1\le k \le n$ (resp., $1\le k \le n-1$ ). The family of $k$-subsets of $G$ (resp., $G^{*}$) that sum up to a given element $x\in G$ is denoted by $ \mathcal{B} _{k}^{x} $ (resp., $ \mathcal{B} _{k}^{x, *} $), and we define $b_{k}^{x}=\left | \mathcal{B} _{k}^{x} \right | $ (resp., $b_{k}^{x,*}=\left | \mathcal{B} _{k}^{x,*} \right | $). For brevity, when $x=0$, we omit the superscript $x$ (for instance, $ \mathcal{B} _{k}^{x} $ is written as $ \mathcal{B} _{k} $).\par
		Additionally, for any $y\in G$ (resp., $G^{*}$) and $x\in G$, the number of $k$-subsets in $ \mathcal{B} _{k}^{x} $ (resp., $ \mathcal{B} _{k}^{x, *} $) that contain $y$ is denoted by $r_{k}^{x}(y)$ (resp., $r_{k}^{x,*}(y)$). When $S$ is a set in $G$ (resp., $G^{*}$), we denote by $r_{k}^{x}(S)$ (resp., $r_{k}^{x,*}(S)$) the number of blocks in $ \mathcal{B} _{k}^{x} $ (resp., $ \mathcal{B} _{k}^{x, *} $) that contain $S$.
	\end{definition}
	The closed forms of $b_{k}^{x}$ and $b_{k}^{x,*}$ for a general finite abelian group $G$ are presented in \cite{K2013}, as stated in the following theorem. First we introduce two notations as follows.\\
	Let exp(G) be the exponent of a group $G$, and for $x\in G$ we define 
	$$e(x)=\max\{d: d|\exp(G), x\in dG\}.$$ 
	For an integer $d$, we denote  the $d$-torsion of $G$ by $$G[d]=\{g\in G:dg=0\}.$$ 
	
	\begin{theorem}\cite{K2013}\label{thm bb}
		Let $G$ be an abelian group of order $n$ and let $\mu$ be the M\"{o}bius function. For $1\le k\le n$, we have:
		$$b_{k}^{x}=\frac{1}{n}\sum_{s|\gcd(\exp(G),k)}(-1)^{k+k/s}\binom{n/s}{k/s}\sum_{d|\gcd(e(x),s)}\mu\left(\frac{s}{d}\right)\#G[d],$$
		and for any $1\le k\le n-1$, we have:
		$$b_{k}^{x,*}=\frac{1}{n}\sum_{s|\exp(G)}(-1)^{k+\lfloor k/s\rfloor}\binom{n/s-1}{\lfloor k/s\rfloor}\sum_{d|\gcd(e(x),s)}\mu\left(\frac{s}{d}\right)\#G[d].$$
		Thus, for $x,y\in G$, if $e(x)=e(y)$, then $b_{k}^{x}=b_{k}^{y}$, $b_{k}^{x,*}=b_{k}^{y,*}$.
	\end{theorem}
	
	When we discuss the incidence structure $(G, \mathcal{B} _{k}^{x})$ as a design, the trivial case of $ \mathcal{B} _{k}^{x}=\emptyset $ must be avoided. The following theorem shows that $ \mathcal{B} _{k}^{x}=\emptyset $ only in the some trivial cases.
	\begin{theorem}\cite{K2013}\label{thm 3c}
		For an abelian group $G$ with $n$ elements and $1\le k\le n-1$, $ \mathcal{B} _{k}^{x}=\emptyset $ if and only if one of the following condition holds:
		\begin{enumerate}
			\item[(i)] \( k = 2, \exp(G) = 2 \) and \( x = 0 \);
			\item[(ii)] \( k = n - 2 \geq 2, \exp(G) = 2 \) and \( x = 0 \);
			\item[(iii)] \( k = n \) and \( x \neq \displaystyle\sum_{g \in G[2]} g \).
		\end{enumerate}
	\end{theorem}
	
	Throughout the paper, for a prime $p$, and any integer $N$, we denote by $\nu_{p}(N)$ the $p$-adic valuation of $N$, defined as $$\nu_{p}(N)=\begin{cases}
		\max\{s\in\mathbb{N}^{+}:p^s\mid N\} & \text{ if } N\ne0, \\
		\infty  & \text{ if } N=0.
	\end{cases}$$
	In addition, we interpret $p^{\infty}$ as $0$, whenever it appears in the proof. 
	\section{Results on Abelian $p$-Groups}\label{sec p}
	Let $p$ be a prime number. In a $p$-group $G$, the order of each element is a power of $p$, which is equivalent to saying that exp($G$)$=p^{m}$, where $m$ is a nonzero integer. In this section, we determine all pairs $(k,x)$ for which $(G, \mathcal{B} _{k}^{x})$ supports a $1$-design, given that exp($G$)$=p^{m}$.\par
	We begin by introducing the following lemma, which is applicable to any finite abelian group.
	\begin{lemma}\label{lem rb}
		Let $G$ be a finite abelian group, $k\ge 2$ and let the parameters $(k,x)$ satisfy $ \mathcal{B} _{k}^{x}\ne\emptyset $. Then $(G, \mathcal{B} _{k}^{x})$ is a $1$-$(n, k, r)$ design if and only for any $y\in G, b_{k-1}^{x-ky,*}=r$, where $r$ is a constant independent of $y$.
	\end{lemma}
	\begin{proof} 
		Consider the map $g\mapsto g-y$, which is a permutation of $G$, and it induces a bijection between the $k$-subsets in $\mathcal{B} _{k}^{x}$ containing $y$ and the $k$-subsets in $\mathcal{B} _{k}^{x-ky}$ containing $0$, hence we have $r_{k}^{x}(y)=b_{k-1}^{x-ky,*}.$ By definition, $(G, \mathcal{B} _{k}^{x})$ is a $1$-$(n, k, r)$ design if and only if for each element $y$,  the number of blocks containing $y$ is the same, this implies $b_{k-1}^{x-ky,*}=r.$  
	\end{proof}
	
	For any abelian $p$-group $G$ whose exponent is a power of the prime $p$, the structure theorem of finite abelian groups implies the following group isomorphism,
	$$G\cong \mathbb{Z}_{p^{t_{1}}} \oplus \mathbb{Z}_{p^{t_{2}}} \oplus\cdots \oplus\mathbb{Z}_{p^{t_{m}}},$$
	where $1\le t_{1}\le\cdots\le t_{m}$ and $\vert G \vert = p^{t_1}\cdot p^{t_2}\cdots p^{t_m}$.\\
	From now on, in this section we fix an abelian $p$-group $G$ of order $n$, and assume it is isomorphic to $ \mathbb{Z}_{p^{t_{1}}} \oplus \mathbb{Z}_{p^{t_{2}}} \oplus\cdots \oplus\mathbb{Z}_{p^{t_{m}}}$, with $1\le t_{1}\le\cdots\le t_{m}$.
	\begin{lemma}\label{lem ne}
		Let $1\le k\le n$ be an integer, for any $g\in G$ satisfying $e(x)=p^{t_{m}-1}$, we have $b_{k}^{*}\ne b_{k}^{g,*}.$
	\end{lemma}
	\begin{proof}
		By Theorem \ref{thm bb}, for any $x\in G$, we have
		$$b_{k}^{x,*}=\frac{1}{n}\sum_{s|p^{t_{m}}}(-1)^{k+\lfloor k/s\rfloor}\binom{n/s-1}{\lfloor k/s\rfloor}\sum_{d|\gcd(e(x),s)}\mu\left(\frac{s}{d}\right)\#G[d].$$
		For the two cases $x=0$ and $x=g$, namely,
		$$b_{k}^{*}=\frac{1}{n}\sum_{s|p^{t_{m}}}(-1)^{k+\lfloor k/s\rfloor}\binom{n/s-1}{\lfloor k/s\rfloor}\sum_{d|\gcd(p^{t_{m}},s)}\mu\left(\frac{s}{d}\right)\#G[d],$$ 
		$$b_{k}^{g,*}=\frac{1}{n}\sum_{s|p^{t_{m}}}(-1)^{k+\lfloor k/s\rfloor}\binom{n/s-1}{\lfloor k/s\rfloor}\sum_{d|\gcd(p^{t_{m}-1},s)}\mu\left(\frac{s}{d}\right)\#G[d],$$
		observe that as $s$ runs through the factors of $p^{t_{m}}$, they only differ in the term $s=p^{t_{m}}$:
		$$\sum_{d|p^{t_{m}}}\mu\left(\frac{p^{t_{m}}}{d}\right)\#G[d]=n+\sum_{d|p^{t_{m}-1}}\mu\left(\frac{p^{t_{m}}}{d}\right)\#G[d].$$
		Thus the lemma is proved by the following equation:
		$$b_{k}^{*}= b_{k}^{g,*}+\frac{1}{n}(-1)^{k+\lfloor k/p^{t_{m}}\rfloor}\binom{n/p^{t_{m}}-1}{\lfloor k/p^{t_{m}}\rfloor}\cdot n=b_{k}^{g,*}+(-1)^{k+\lfloor k/p^{t_{m}}\rfloor}\binom{n/p^{t_{m}}-1}{\lfloor k/p^{t_{m}}\rfloor}.$$
	\end{proof}
	\begin{theorem}\label{thm p}
		Let $p$ be an odd prime and let $G$ be an abelian $p$-group isomorphic to $ \mathbb{Z}_{p^{t_{1}}} \oplus \mathbb{Z}_{p^{t_{2}}} \oplus\cdots \oplus\mathbb{Z}_{p^{t_{m}}}$, where $1\le t_{1}\le \cdots\le t_{m}$. For any $1\le k\le n$, and $x=(x_{1}, x_{2},\cdots, x_{m})\in G$, with $x_{i}\in \mathbb{Z}_{p^{t_{i}}}  $, $(G, \mathcal{B} _{k}^{x})$ is a $1$-$(n, k, r)$ design if and only if $p\mid k$ and the pair $(k,x)$ satisfies one of the following conditions:
		\begin{enumerate}
			\item[(i)] $p^{t_{m}}|k$, $k\ne n$, and $x\in G$ is an arbitrary element, or $k=n$, $x=0.$
			\item[(ii)] $k\ne n$, and there exist at least one $ i, 1\le i\le m$, such that $p\nmid x_{i}.$ 
			\item[(iii)] $k\ne n$, $p \mid x_{i}$ for all $1\le i\le m$, and $ \max\left \{ \nu_{p}(k)-\nu_{p}^{i}(x_{i})\mid 1\le i\le m \right \} \ge1$, where $\nu_{p}^{i}$ is the $p$-adic valuation restricted to $\mathbb{Z}_{p^{t_{i}}}$ defined as
			$$\nu_{p}^{i}(x)=
			\begin{cases}
				\nu_{p}(x) & \text{ if } x\ne0\pmod {p^{t_{i}}},\\
				\infty  & \text{ if } x=0\pmod {p^{t_{i}}}.
			\end{cases}$$
		\end{enumerate}
	\end{theorem}
	\begin{proof}
		Since it is evident that $(G, \mathcal{B} _{k}^{x})$ cannot support a $1$-design when $k=1$, so we assume $k\ge 2$. The order of $G$ is odd, as $p$ is an odd prime. We have $$\sum_{g \in G[2]}g=\sum_{g \in G}g=0.$$
		Hence, by Theorem \ref{thm 3c}, $ \mathcal{B} _{k}^{x}=\emptyset $ only when $k=n$ and $x\ne 0$.\\ 
		In the case $ \mathcal{B} _{k}^{x}\ne\emptyset $, Lemma \ref{lem rb} implies that for any $y=(y_{1}, y_{2},\cdots, y_{m})\in G$, $r_{k}^{x}(y)=b_{k-1}^{x-ky,*}$, and $(G, \mathcal{B} _{k}^{x})$ is a $1$-$(n, k, r)$ design if and only if $b_{k-1}^{x-ky,*}=r$ is a constant independent of $y$. Now we consider the following two cases:\\
		When $k=n$, then $(G, \mathcal{B} _{k}^{x})$ is a $1$-design if and only if $(k,x)=(n,0)$. When $k\ne n$ and exp$(G)=p^{t_{m}}\mid k$, then $b_{k-1}^{x-ky,*}=b_{k-1}^{x,*}$, which is a constant independent of $y$, thus $(G, \mathcal{B} _{k}^{x})$ is a $1$-design.\\ 
		We claim that $(G, \mathcal{B} _{k}^{x})$ is a $1$-design only if $p \mid  k$. If possible assume that $p \nmid  k$, then as $y_{i}$ runs through $\mathbb{Z}^{p_{t_{i}}}$, $x_{i}-ky_{i}$ also runs through $\mathbb{Z}^{p_{t_{i}}}$, thus $x-ky$ runs through $G$ when $y$ runs through $G$, that is,
		$$G=\left \{(x_{1}-ky_{1},\cdots,x_{m}-ky_{m}) \mid (y_{1},\cdots,y_{m})\in G \right \}.$$
		Therefore, there exist $y$ and $y'$ such that $x-ky=0$ and $x-ky'=(0,\cdots,0, p^{t_{m}-1} )$, while by Lemma \ref{lem ne} $b_{k-1}^{x-ky,*}\ne b_{k-1}^{x-ky',*} $, hence $(G, \mathcal{B} _{k}^{x})$ is not a $1$-design when $p \nmid  k$.\\ 
		Thus, in the remainder of the proof, we assume $p\mid k$. Given $p\mid k$, if there exists a coordinate $x_{i}$ of $x$ such that $p\nmid x_{i}$, then $p\nmid x_{i}-ky$ for any $y\in G$. It follows that $e(x-ky)=1$ for any $y\in G$. By  Theorem \ref{thm bb}, $b_{k-1}^{x-ky,*}$ is a constant as $y$ runs through $G$, then in this case $(G, \mathcal{B} _{k}^{x})$ is a $1$-design. If $p\mid k$, $p^{t_{m}}\nmid k$ and $p\mid x_{i}$ for any $1\le i\le m$, then if $ \max\left \{ \nu_{p}(k)-\nu_{p}^{i}(x_{i})\mid 1\le i\le m \right \} \le0$, that is, $\nu_{p}(k)\le\nu_{p}^{i}(x_{i})$ for all $1\le i \le m$, in this case, we have
		$$x_{i}-ky_{i}=w_{i_{1}}p^{\nu_{p}^{i}(x_{i})}-w_{i_{2}}p^{\nu_{p}(k)}y_{i}=p^{\nu_{p}(k)}(w_{i_{1}}p^{\nu_{p}^{i}(x_{i})-\nu_{p}(k)}-w_{i_{2}}y_{i}),$$
		where $\gcd(w_{i_{j}},p)=1$, for $j=1,2$, and we further notice that when $y_{i}$ runs through $\mathbb{Z}_{p^{t_{i}}}$, $w_{i_{2}}y_{i}$ also runs through $\mathbb{Z}_{p^{t_{i}}}$, hence,
		$$\mathbb{Z}_{p^{t_{i}}}=\left \{ w_{i_{1}}p^{\nu_{p}^{i}(x_{i})-\nu_{p}(k)}-w_{i_{2}}y_{i}\mid y_{i}\in \mathbb{Z}_{p^{t_{i}}}   \right \}.$$
		Therefore, there exist $y,y'\in G$, such that $x-ky=0$ and $e(x-ky')=p^{t_{m}-1}$, which implies that $(G, \mathcal{B} _{k}^{x})$ is not a $1$-design. On the contrary, if $p\mid k$, $p^{t_{m}}\nmid k$ and $p\mid x_{i}$ for any $1\le i\le m$, and $ \max\left \{ \nu_{p}(k)-\nu_{p}^{i}(x_{i})\mid 1\le i\le m \right \} \ge1$, then without loss of generality, let $\nu_{p}(k)-\nu_{p}^{i}(x_{i})\ge1$, $1\le i\le j$ and 
		$\nu_{p}(k)\le\nu_{p}^{i}(x_{i})$, $j+1\le i\le m$. In this case, we have
		$$x_{i}-ky_{i}=
		\begin{cases}
			p^{\nu_{p}^{i}(x_{i})}(w_{i_{1}}-w_{i_{2}}p^{\nu_{p}(k)-\nu_{p}^{i}(x_{i})}y_{i})	& \text{ if } 1\le i\le j, \\
			p^{\nu_{p}(k)}(w_{i_{1}}p^{\nu_{p}^{i}(x_{i})-\nu_{p}(k)}-w_{i_{2}}y_{i})	& \text{ if } j+1\le i\le m,
		\end{cases}$$
		with $\gcd(w_{i_{j}},p)=1$, for $1\le i\le m$ and $j=1,2$. Therefore, for any $y\in G$, 
		$$e(x-ky)=p^{\min \{\nu_{p}(x_{i}-ky_{i})\mid 1\le i\le j  \}}=p^{\min \{\nu_{p}^{i}(x_{i})\mid 1\le i\le j \}}.$$
		Hence, in this case $(G, \mathcal{B} _{k}^{x})$ is a $1$-design.
		Combine these cases together, $(G, \mathcal{B} _{k}^{x})$ is a $1$-design if and only if $p\mid k$ and one of the three conditions holds.
	\end{proof}
	For an odd prime $p$, we have determined all pairs $(k,x)$ such that $(G, \mathcal{B} _{k}^{x})$ supports a $1$-design over the abelian $p$-group $ \mathbb{Z}_{p^{t_{1}}} \oplus \mathbb{Z}_{p^{t_{2}}} \oplus\cdots \oplus\mathbb{Z}_{p^{t_{m}}}$ in Theorem \ref{thm p}. The following result addresses the case of $p=2$.
	\begin{proposition}\label{co p=2}
		For an abelian $2$-group $G \cong\mathbb{Z}_{2^{t_{1}}} \oplus \mathbb{Z}_{2^{t_{2}}} \oplus\cdots \oplus\mathbb{Z}_{2^{t_{m}}}$, with $1\le t_{1}\le\cdots\le t_{m}$, $1\le k\le n$, and $x=(x_{1}, x_{2},\cdots, x_{m})\in G$. When $ t_{m}>1$, $(G, \mathcal{B} _{k}^{x})$ is a $1$-$(n, k, r)$ design if and only if the pair of parameters $(k,x)$ satisfies one of the following conditions:
		\begin{enumerate}
			\item[(i)] $2^{t_{m}}|k$, $k\ne n$, and $x\in G$ is an arbitrary element, or $k=n$, $x=\begin{cases}
				\frac{n}{2}  & \text{ if } m=1, \\
				0 & \text{ if } m>1.
			\end{cases}$
			\item[(ii)] $k\ne n$, and there exist at least one $ i, 1\le i\le m$, such that $2\nmid x_{i}.$ 
			\item[(iii)] $k\ne n$, $2 \mid x_{i}$ for any $1\le i\le m$, and $k$ satisfies $ \max\left \{ \nu_{2}(k)-\nu_{2}^{i}(x_{i})\mid 1\le i\le m \right \} \ge1$, where $\nu_{2}^{i}$ is a function over $\mathbb{Z}_{2^{t_{i}}}$ defined as
			$$\nu_{2}^{i}(x)=
			\begin{cases}
				\nu_{2}(x)& \text{ if } x\ne0, \\
				2^{t_{i}}& \text{ if } x=0.
			\end{cases}$$
		\end{enumerate}
		When $t_{m}=1$, $m>1$, $(G, \mathcal{B} _{k}^{x})$ is a $1$-$(n, k, r)$ design if and only if the pair $(k,x)$ belongs to the set $\Omega_{2}\setminus \left \{ (2,0),(n-2,0) \right \} $, where $\Omega_{2}$ is the set determined by conditions $(i)$, $(ii)$ and $(iii)$. When $t_{m}=1$, $m=1$, and $G=\mathbb{Z}_{2}$, the case is trivial.
	\end{proposition}
	\begin{proof} When $p=2$, we have 
		$$\sum_{g \in G[2]}g=\sum_{g \in G}g=\begin{cases}
			\frac{n}{2}  & \text{ if } m=1,\\
			0 & \text{ if } m>1.
		\end{cases}$$
		By Theorem \ref{thm 3c}, $ \mathcal{B} _{k}^{x}=\emptyset$ if and only if $(k,x)$ is $(2,0)$, $(n-2,0)$ or $(n, x)$ with $x\ne \sum_{g \in G[2]}g $, hence the proof follows the same approach as Theorem \ref{thm p}, where the only additional consideration is the avoidance of the case $ \mathcal{B} _{k}^{x}=\emptyset$.
	\end{proof}
	It is important to note that the result presented in \cite{P2023} can be directly derived from Theorem \ref{thm p} in the following result.
	\begin{corollary}\label{coro in}
		Let $G=\mathbb{F}_{p}^{m}$, where $p$ is an odd prime. Then the incidence structure $(G,\mathcal{B} _{k}^{x} )$ is a $1$-design if and only if $p\mid k$, $k\ne p^{m}$, $x\in G$ is an arbitrary element, or $k=p^{m}$, $x=0$.
	\end{corollary}
	\begin{proof}
		In the case $G=\mathbb{F}_{p}^{m}=\underbrace{\mathbb{Z}_{p}\oplus\cdots\oplus\mathbb{Z}_{p}}_{\text{$m$ times}}$, the condition $p\mid k$ is necessary and by $(i)$ in Theorem \ref{thm p} where $x\in G$ can be any element except in the trivial case where $k=p^{m}$, in which $x$ must be $0$.
	\end{proof}
	\begin{remark}\label{rem}
		
		In many problems related to subset sums in a group, the sum $x$ is often fixed to be $0$, for example, the well-known \emph{zero-sum problems} \cite{N1996, T2006}. It is interesting that when $x=0$, by Theorem \ref{thm p}, the condition for $(G,\mathcal{B} _{k} )$ to support a $1$-design is concise: $p^{t_{m}}\mid k$, that is, the exponent $\exp(G)$ divides $k$. Further, one may ask when $(G,\mathcal{B} _{k}^{x} )$ supports a $2$-design. In \cite{P2023} the author proved that $(G,\mathcal{B} _{k}^{x})$ supports a $2$-design if and only if $p\mid k$ and $x=0$, providing $G=\mathbb{F}_{p}^{m}$. Although the proof for general case in \cite{P2023} is highly technical, but in the case $x=0$, alternatively, there is an elegant proof, thanks to the structure of finite fields, as shown in the following proposition.
	\end{remark}
	\begin{proposition}\cite{P2023}\label{x=0}
		Let $G=\mathbb{F}_{p}^{m}=\underbrace{\mathbb{Z}_{p}\oplus\cdots\oplus\mathbb{Z}_{p}}_{\text{$m$ times}}$, when $x=0$ and $\mathcal{B} _{k} \ne \emptyset$, $(G,\mathcal{B} _{k} )$ is a $2$-design if and only if $p\mid k$.
	\end{proposition}
	\begin{proof}
		It remains to prove the condition $p\mid k$ is sufficient, in this case, note that any affine mapping $L(x)+b$ is a permutation of the blocks in $\mathcal{B} _{k}$, where $L$ is an invertible linear map, and $b$ is a constant vector. Note that the group of affine maps $\mathrm{AGL}(m,p)$ acts $2$-transitively on $\mathbb{F}_{p}^{m}$, that is, for any two distinct pairs $(x_{1}, y_{1})$ and $(x_{2}, y_{2})$, there exists an affine map $L(x)+b$, such that $L(x_{1}, y_{1})=(x_{2}, y_{2})$, then the affine map $L(x)+b$ induces a one to one correspondence of blocks contain $(x_{1}, y_{1})$ and blocks contain $(x_{2}, y_{2})$, which shows $(G,\mathcal{B} _{k} )$ supports a $2$-design.
	\end{proof}
	While the condition for $(G,\mathcal{B}_{k}^{x} )$ to support a $2$-design is fully characterized when $G=\mathbb{F}_{p}^{m}$, the extension to non-elementary abelian $p$-groups may fail. For any non-elementary abelian $p$-group, that is, $G=\mathbb{Z}_{p^{t_{1}}} \oplus \mathbb{Z}_{p^{t_{2}}} \oplus\cdots \oplus\mathbb{Z}_{p^{t_{m}}}$ with some $t_{i}>1$, computational verification of such groups via MAGMA suggests that the incidence structure $(G,\mathcal{B}_{k}^{x} )$ probably fails to be a $2$-design for any pair $(k,x)$, except the trivial case of $k=|G|$. This failure may be fundamentally attributed to the collapse of certain Symmetric hierarchy in non-elementary $p$-abelian groups. Therefore, we propose the following conjecture.
	\begin{conjecture}\label{conjecture}
		Let $G=\mathbb{Z}_{p^{t_{1}}} \oplus \mathbb{Z}_{p^{t_{2}}} \oplus\cdots \oplus\mathbb{Z}_{p^{t_{m}}}$ be a non-elementary $p$-abelian group, with some $t_{i}>1$. For any $(k,x)$, where $1\le k\le n-1$ and $x\in G$, the incidence structure $(G,\mathcal{B}_{k}^{x} )$ is not a $2$-design.
	\end{conjecture}
	
	\section{Results on Abelian Groups of Exponent $pq$}\label{sec pq}
	This section determines the conditions for $(G, \mathcal{B} _{k}^{x})$ to be a $1$-$(n, k, r)$ design over any group with exponent $pq$, where $p$ and $q$ are two primes and $p<q$.\\ 
	The structure of finite abelian groups implies that an abelian group of exponent $pq$ is isomorphic to either $$\underbrace{\mathbb{Z}_{p}\oplus\cdots\oplus\mathbb{Z}_{p}}_{\text{$s$ times}}\oplus\underbrace{\mathbb{Z}_{pq}\oplus\cdots\oplus \mathbb{Z}_{pq}}_{\text{$t-s$ times}},$$  or
	$$\underbrace{\mathbb{Z}_{q}\oplus\cdots\oplus\mathbb{Z}_{q}}_{\text{$s$ times}}\oplus\underbrace{\mathbb{Z}_{pq}\oplus\cdots\oplus \mathbb{Z}_{pq}}_{\text{$t-s$ times}},$$ where $0\le s\le t, t\ge 1$. Let $g\in G$, then $$e(g)=\max\{d:d\mid pq,g\in dG\}.$$
	\begin{lemma}\label{lem pq}
		Consider the group $G\cong\underbrace{\mathbb{Z}_{p}\oplus\cdots\oplus\mathbb{Z}_{p}}_{\text{$s$ times}}\oplus\underbrace{\mathbb{Z}_{pq}\oplus\cdots\oplus \mathbb{Z}_{pq}}_{\text{$t-s$ times}}$. Let $k\ge2$ and let $g_{1}, g_{2}, g_{3}\in G$. Then
		\begin{enumerate}
			\item [(i)] For a fixed $k\ge 2$, $\#\left \{ b_{k-1}^{*,g}\mid g\in G \right \}\ge2$.
			\item [(ii)] If $e(g_{1})=q$, $e(g_{2})=1$, then $b_{k-1}^{g_{1},*}= b_{k-1}^{g_{2},*}$ if and only if $k\le q$ or $k\ge n-q+1$.
			\item [(iii)] If $e(g_{3})=p$, $e(g_{2})=1$, then $b_{k-1}^{g_{3},*}= b_{k-1}^{g_{2},*}$ if and only if $k\le p$ or $k\ge n-p+1$.
		\end{enumerate}
	\end{lemma}
	\begin{proof}
		We pick four elements $g_{1}, g_{2}, g_{3},g_{4}$ from $G$, where $e(g_{1})=q$, $e(g_{2})=1$, $e(g_{3})=p$, $e(g_{4})=pq$ ($g_{4}=0$). We then prove $(i)$ by showing that if $b_{k-1}^{g_{1},*}= b_{k-1}^{g_{2},*}$, then $b_{k-1}^{g_{3},*}\ne b_{k-1}^{g_{4},*}$. By Theorem \ref{thm bb}, $b_{k-1}^{g_{1},*}= b_{k-1}^{g_{2},*}$ if and only if 
		\begin{align*}
			&\frac{1}{n} [\binom{n-1}{k-1}+(-1)^{k-1+\left \lfloor (k-1) /p\right \rfloor } \binom{n/p-1}{\left \lfloor (k-1)/p \right \rfloor } (-1)+(-1)^{k-1+\left \lfloor (k-1) /q\right \rfloor } \\&\binom{n/q-1}{\left \lfloor (k-1)/q \right \rfloor } (q^{t}-1) +(-1)^{k-1+\left \lfloor (k-1) /pq\right \rfloor }\binom{n/pq-1}{\left \lfloor (k-1)/pq \right \rfloor }(1-q^{t}) ]\\
			=&
			\frac{1}{n} [\binom{n-1}{k-1}+(-1)^{k-1+\left \lfloor (k-1) /p\right \rfloor } \binom{n/p-1}{\left \lfloor (k-1)/p \right \rfloor } (-1)+(-1)^{k-1+\left \lfloor (k-1) /q\right \rfloor } \\&\binom{n/q-1}{\left \lfloor (k-1)/q \right \rfloor } (-1) +(-1)^{k-1+\left \lfloor (k-1) /pq\right \rfloor }\binom{n/pq-1}{\left \lfloor (k-1)/pq \right \rfloor }(1-q^{t}) ].
		\end{align*}
		That is,
		\begin{equation}\label{eq12}
			(-1)^{\left \lfloor (k-1) /q\right \rfloor }\binom{n/q-1}{\left \lfloor (k-1)/q \right \rfloor }=(-1)^{\left \lfloor (k-1) /pq\right \rfloor}\binom{n/pq-1}{\left \lfloor (k-1)/pq \right \rfloor }.
		\end{equation}
		Similarly, $b_{k-1}^{g_{3},*}=b_{k-1}^{g_{4},*}$ if and only if 
		\begin{equation}\label{eq34}
			(-1)^{\left \lfloor (k-1) /q\right \rfloor }\binom{n/q-1}{\left \lfloor (k-1)/q \right \rfloor }=(-1)^{\left \lfloor (k-1) /pq\right \rfloor}\binom{n/pq-1}{\left \lfloor (k-1)/pq \right \rfloor }(1-p^{t+s}).
		\end{equation}
		However, it is impossible for Eqs. (\ref{eq12}) and (\ref{eq34}) to hold simultaneously, which implies that $$\#\left \{ b_{k-1}^{*,g}\mid g\in G \right \}\ge2.$$
		For $(ii)$ and $(iii)$, we only prove $(ii)$, then $(iii)$ is obtained by replacing $q$ by $p$. If $e(g_{1})=q$, $e(g_{2})=1$, then 
		it is obvious that Eq. (\ref{eq12}) holds when $k\le q$ or $k\ge n-q+1$, then we denote 
		$$f(n,k)=\frac{\binom{n/q-1}{\left \lfloor (k-1)/q \right \rfloor }}{\binom{n/pq-1}{\left \lfloor (k-1)/pq \right \rfloor }}.$$
		Subsequently, it is evident that $f(n,k)=f(n,n+1-k)$ through straightforward calculations. Thus, it suffices to prove that when $q+1\le k\le \frac{n+1}{2}$, $f(n,k)>1$, so that $b_{k-1}^{g_{1},*}\ne b_{k-1}^{g_{2},*}$. When $q+1\le k\le \frac{n+1}{2}$, then $1\le\left \lfloor \frac{k-1}{q} \right \rfloor\le \left \lfloor \frac{n-1}{2q}\right\rfloor <\frac{1}{2}\left \lfloor\frac{n-1}{q}\right \rfloor-1<\frac{1}{2}(\frac{n}{q}-1)$, and we have
		\begin{equation}\label{eqAB}
			f(n,k)=\underbrace{\frac{\binom{n/q-1}{\left \lfloor (k-1)/q \right \rfloor } }{\binom{n/q-1}{\left \lfloor (k-1)/pq \right \rfloor }}}_{A}\cdot\underbrace{ \frac{\binom{n/q-1}{\left \lfloor (k-1)/pq \right \rfloor } }{\binom{n/pq-1}{\left \lfloor (k-1)/pq \right \rfloor }}}_{B}.
		\end{equation}
		In Eq. (\ref{eqAB}) $B\ge1$ is obvious, then by $\left \lfloor\frac{k-1}{pq} \right \rfloor\le \frac{1}{p}\left \lfloor\frac{k-1}{q} \right \rfloor$ and $1\le\left \lfloor\frac{k-1}{q} \right \rfloor<\frac{1}{2}(\frac{n}{q}-1)$, $A>1$. Hence, $f(n,k)>1$, then $b_{k-1}^{g_{1},*}\ne b_{k-1}^{g_{2},*}$ when $q+1\le k \le n-q$, which completes the proof.
	\end{proof}
	\begin{theorem}\label{thm pq}
		Let $G\cong \underbrace{\mathbb{Z}_{p}\oplus\cdots\oplus\mathbb{Z}_{p}}_{\text{$s$ times}}\oplus\underbrace{\mathbb{Z}_{pq}\oplus\cdots\oplus \mathbb{Z}_{pq}}_{\text{$t-s$ times}}$, $|G|=n$, $x=(x_{1}, x_{2},\cdots, x_{t})\in G$. Then $(G, \mathcal{B}_{k}^{x})$ is a $1$-$(n, k, r)$ design if and only if one of the following conditions holds:
		\begin{enumerate}
			\item[(i)]  $k=n$, $x=0.$
			\item[(ii)] $k\ne n$, $pq\mid k$, and $x$ is an arbitrary element in $G$.
			\item[(iii)] $p\mid k$, there exists at least one $x_{i}$ such that $p\nmid x_{i}$, and $k$ satisfies $k\le q-1$ or $k\ge n-q+1$.
			\item[(iv)] $p\mid x_{i}$ for any $1\le i\le t$, and $k$ satisfies $p\mid k$, $q\nmid k$, and
			$$(-1)^{\left \lfloor \frac{k-1}{q} \right \rfloor-\left \lfloor \frac{k-1}{pq} \right \rfloor}\frac{\binom{n/q-1}{\left \lfloor (k-1)/q \right \rfloor }}{\binom{n/pq-1}{\left \lfloor (k-1)/pq \right \rfloor }}=(1-p^{t+s}).$$
			\item[(v)] $q\mid x_{i}$ for any $1\le i\le t$, and $k$ satisfies $q\mid k$, $p\nmid k$, and
			$$(-1)^{\left \lfloor \frac{k-1}{p} \right \rfloor-\left \lfloor \frac{k-1}{pq} \right \rfloor}\frac{\binom{n/p-1}{\left \lfloor (k-1)/p \right \rfloor }}{\binom{n/pq-1}{\left \lfloor (k-1)/pq \right \rfloor }}=(1-q^{t}).$$
		\end{enumerate}
	\end{theorem}
	\begin{proof}
		By Theorem \ref{thm 3c}, $ \mathcal{B} _{k}^{x}=\emptyset$ if and only if $k=n$ and $x\ne 0$. When $k=1$, $ \mathcal{B} _{k}^{x}\ne \emptyset$ is not a $1$-design, so we assume $2\le k<n$, then $ \mathcal{B} _{k}^{x}\ne \emptyset$, and by Lemma \ref{lem rb}, $(G,\mathcal{B} _{k}^{x})$ is a $1$-design if and only if $b_{k-1}^{*,x-ky}$ is a constant for any $y=(y_{1},\cdots, y_{t})\in G$. If $p\nmid k$, $q\nmid k$, then $(G,\mathcal{B} _{k}^{x})$ is not a $1$-design because when $y_{i}$ runs through $G$, $x_{i}-ky_{i}$ also runs through $\mathbb{Z}_{p}$ and when $1\le i\le s$, and runs through $\mathbb{Z}_{pq}$ and when $s\le i\le t$, which follows that $x-ky$ runs through $G$, that is, $G=\{y\in G\mid x-ky\}$, but by Lemma \ref{lem pq} $\#\left \{b_{k-1}^{*,g}\mid g\in G \right \}\ge2$, hence $b_{k-1}^{*,x-ky}$ is not a constant as $y$ runs through $G$. Thus, $(G,\mathcal{B} _{k}^{x})$ is a $1$-design only if $p\mid k$ or $q\mid k$, and if $pq\mid k$, then $x-ky=x$, $b_{k-1}^{*,x-ky}$ is a constant. If $p\mid k$, $q\nmid k$, for $x=(x_{1}, x_{2},\cdots, x_{t})\in G$ we have $$\{\gcd(x_{i}-ky_{i}, pq )\mid y_{i}\in \mathbb{Z}\}=
		\begin{cases}
			\{q,1\}& \text{ if } \nu_{p}(x_{i}) =0,  \\
			\{p,pq\}& \text{ if } \nu_{p}(x_{i}) \ge1.
		\end{cases}$$ 
		Therefore, $$\{e(x-ky)\mid y\in G\}=
		\begin{cases}
			\{q,1\}& \text{ if } \,\exists i, \nu_{p}(x_{i}) =0,  \\
			\{p,0\}& \text{ if } \nu_{p}(x_{i}) \ge1, 1\le i\le t.
		\end{cases}$$
		Then $(iii)$ and $(iv)$ are immediately obtained by Lemma \ref{lem pq}. In the case that $p\nmid k$, $q\mid k$, the proof is almost the same as the case that $p\mid k$, $q\nmid k$, but here $(G,\mathcal{B} _{k}^{x})$ is a $1$-design only if $q\mid x_{i}$ for any $1\le i\le t$ because if $k\le p-1$ or $k\ge n-p+1$, then $q\nmid k$ for $p<q$, hence in this case we only have $(v)$.
	\end{proof}
	In the case $G\cong\underbrace{\mathbb{Z}_{q}\oplus\cdots\oplus\mathbb{Z}_{q}}_{\text{$s$ times}}\oplus\underbrace{\mathbb{Z}_{pq}\oplus\cdots\oplus \mathbb{Z}_{pq}}_{\text{$t-s$ times}}$, the conditions for $(G,\mathcal{B} _{k}^{x})$ to be a $1$-design is exactly the same as the conditions in Theorem \ref{thm pq}, with the same proof of the case $G\cong\underbrace{\mathbb{Z}_{p}\oplus\cdots\oplus\mathbb{Z}_{p}}_{\text{$s$ times}}\oplus\underbrace{\mathbb{Z}_{pq}\oplus\cdots\oplus \mathbb{Z}_{pq}}_{\text{$t-s$ times}}$.
	\section{Results on General Finite Abelian Groups of Cyclic and Non-Cyclic Type} \label{sec cyc}
	In this section, as further exploration of the $t$-designs hold in general finite abelian groups, we propose some interesting results of incidence structure $(G, \mathcal{B} _{k}^{x})$, where $G$ is a cyclic group $\mathbb{Z}_{n}$ and general finite abelian group $\mathbb{Z}_{n_{1}}\oplus\cdots\oplus\mathbb{Z}_{n_{m}}$, which are closely related elliptic curve codes in Section \ref{elliptic curve codes}.
	\subsection{Cyclic Groups }
	In the following proposition, we give a necessary condition for $(\mathbb{Z}_{n}, \mathcal{B} _{k}^{x})$ to be a $1$-design, where the largest prime factor of $n$ is relatively large.
	\begin{proposition}\label{prop cyc1}
		Let $G=\mathbb{Z}_{n}$, $x\in G$, and let $q$ be the largest prime factor of $n$, with $\nu_{q}(n)=t$. If $n\le q^{2t}-1$, then $(G, \mathcal{B} _{k}^{x})$ is not a $1$-design if $\gcd(k,n)= 1$.
	\end{proposition}
	\begin{proof}
		When $\gcd(k,n)= 1$, by Theorem \ref{thm bb}, the number of blocks is $b_{k}^{x}=\binom{n}{k}/n$, if $(G, \mathcal{B} _{k}^{x})$ supports a $1$-$(n,k,r)$ design, then we have the double counting argument: 
		$$n\cdot r=b_{k}^{x}\cdot k=\frac{\binom{n}{k}}{n}\cdot k,$$
		which leads to $n^{2}\mid \binom{n}{k}$, but in the following we will prove it is impossible when $\gcd(k,n)= 1$: with $\nu_{q}(n^{2})=2t$, $n\le q^{2t}-1$, by Legendre's Theorem we have
		\begin{align*}
			\nu_{q}(\binom{n}{k} )=&\nu_{q}(\frac{n!}{k!(n-k)!} )\\=&\sum_{i=1}^{\infty } \left \lfloor \frac{n}{q^{i}}  \right \rfloor-\sum_{i=1}^{\infty }\left \lfloor \frac{k}{q^{i}}  \right \rfloor-\sum_{i=1}^{\infty }\left \lfloor \frac{n-k}{q^{i}}  \right \rfloor\\=&\sum_{i=1}^{2t-1 }\left \lfloor \frac{n}{q^{i}}  \right \rfloor-\left \lfloor \frac{k}{q^{i}}  \right \rfloor-\left \lfloor \frac{n-k}{q^{i}}  \right \rfloor\\\le& \,2t-1\\<&\nu_{q}(n^{2}).
		\end{align*}
		Therefore, $n^{2}\nmid \binom{n}{k}$, and then $(G, \mathcal{B} _{k}^{x})$ is not a $1$-design when $\gcd(n,k)=1$, which completes the proof.
	\end{proof}
	
	In the subsequent result, we examine the two incidence structures $(G, \mathcal{B} _{k}^{x_{1}})$ and $(G, \mathcal{B} _{k}^{x_{2}})$, as well as the connections of designs that exist in them.
	\begin{theorem}\label{thm cyc2}
		Let $G=\mathbb{Z}_{n}$, $2\le k\le n-1$, and $x_{1}, x_{2}\in G$. If $\gcd (x_{1},n)=\gcd(x_{2},n)$ then $(G, \mathcal{B} _{k}^{x_{1}})$ and $(G, \mathcal{B} _{k}^{x_{2}})$ are either both $1$-designs or neither is a $1$-design.
	\end{theorem}
	\begin{proof}
		Assume that $n=\displaystyle\prod_{k=1}^{t} p_{k}^{n_{k}}$, where $p_{k}$ are distinct primes, $n_{k}\ge 1$, $1\le k\le t$. Moreover, let $\gcd(x_{1},n)=\gcd(x_{2},n)=\displaystyle\prod_{l=1}^{m}p_{i_{l}}^{\overline{n_{i_{l}}}}$, where $1\le i_{1}<i_{2}<\cdots<i_{m}\le t$, and $\overline{n_{i_{l}}}\le n_{i_{l}}$, $1\le l\le m$. We will first change the form of $x_{j}$, $j=1,2$, before it, note that
		\begin{equation}\label{eq 2j}
			\begin{cases}
				\nu_{p_{i_{l}}}(\frac{x_{j}}{\gcd(x_{j},n)}) \ge0 & \text{ if } \overline{n_{i_{l}}}={n_{i_{l}}},\\
				\nu_{p_{i_{l}}}(\frac{x_{j}}{\gcd(x_{j},n)})=0 & \text{ if }\overline{n_{i_{l}}}<{n_{i_{l}}},
			\end{cases}
		\end{equation}
		where $j=1,2$ and $1\le l\le m$, then, for $j=1,2$:
		\begin{align*}
			x_{j}=&\frac{x_{j}}{\gcd(x_{j},n)}\gcd(x_{j},n)-n\\=&\frac{x_{j}}{\gcd(x_{j},n)}\displaystyle\prod_{l=1}^{m}p_{i_{l}}^{\overline{n_{i_{l}}}}-\displaystyle\prod_{k=1}^{t} p_{k}^{n_{k}}\\=&\displaystyle\prod_{l=1}^{m}p_{i_{l}}^{\overline{n_{i_{l}}}}(\frac{x_{j}}{\gcd(x_{j},n)}\gcd(x_{j},n)-\displaystyle\prod_{l=1}^{m}p_{i_{l}}^{n_{i_{l}}-\overline{n_{i_{l}}}}\cdot\displaystyle\prod_{k\ne i_{l}}p_{k})\\\triangleq&\displaystyle\prod_{l=1}^{m}p_{i_{l}}^{\overline{n_{i_{l}}}}w_{j}\\=&\gcd(x_{j},n)w_{j}.
		\end{align*}
		By Eq. (\ref{eq 2j}), $\gcd(w_{j},n)=1$. Therefore, for any $1\le k\le n$ and any $y\in G$, we have
		\begin{equation}\label{eq run}
			x_{2}-k\frac{w_{2}}{w_{1}}y=\frac{w_{2}}{w_{1}}(x-ky),
		\end{equation}
		where $\gcd(\frac{w_{2}}{w_{1}},n)=1$, this is because $\gcd(w_{1},n)=1$. By B\'ezout's Theorem, there are $w_{1}^{-1}, t\in\mathbb{Z}_{n} $ such that $w_{1}w_{1}^{-1}+nt=1$, then $\gcd(w_{1}^{-1},n)=1$, which follows that $\gcd(\frac{w_{2}}{w_{1}},n)=1$. Note that in the group $G=\mathbb{Z}_{n}$, except the trivial case $n=2$, for $x\in G$ we have $e(x)=\gcd(x,n)$ and 
		\begin{equation}\label{eq cyc b}
			b_{k}^{x,*}=\frac{1}{n}\sum_{s|n}(-1)^{k+\lfloor k/s\rfloor}\binom{n/s-1}{\lfloor k/s\rfloor}\sum_{d|\gcd(e(x),s)}\mu\left(\frac{s}{d}\right)d.
		\end{equation}
		Then by Eq. (\ref{eq run}), $$\{e(x_{1}-ky)\mid y\in G \}=\{e(x_{2}-ky)\mid y\in G \}.$$
		Therefore, for $2\le k\le n-1$, $b_{k-1}^{x_{1}-ky,*}$ and $b_{k-1}^{x_{2}-ky,*}$ are either both are constants or neither is a constant. By Lemma \ref{lem rb} $(G, \mathcal{B} _{k}^{x_{1}})$ and $(G, \mathcal{B} _{k}^{x_{2}})$ are either both $1$-designs or neither is a $1$-design.
	\end{proof}
	\subsection{Non-Cyclic Abelian Groups}
	This section considers general abelian group $G$ of non-cyclic type. We will characterize the connection between $t$-designs hold in $(G,\mathcal{B}_{k}^{x} )$ and $(t-1)$-designs hold in $(G^{*},\mathcal{B}_{k}^{x,*} )$, which applies to the case when designs on the structure $(G^{*},\mathcal{B}_{k}^{x,*} )$ is considered. For example, in Section \ref{elliptic curve codes}, we will deal with designs on $(E(\mathbb{F}_{q})^{*},\mathcal{B}_{k}^{x,*} )$, where $E(\mathbb{F}_{q})$ is the group of rational points of an elliptic curve $E$.
	\begin{theorem}\label{t and t-1}
		Suppose $G$ is a finite abelian group with exponent $\exp G$. Let $\exp G$ divides $k$, and let $k<n$, $x\in G$. Then we have
		\begin{enumerate}
			\item [(i)] If $(G,\mathcal{B}_{k}^{x} )$ is a $t$-design ($t\ge 2$), then $(G^{*}, \mathcal{B}_{k}^{x,*} )$ is a $(t-1)$-design.
			\item [(ii)] Conversely, If $(G^{*},\mathcal{B}_{k}^{x,*} )$ and $(G,\mathcal{B}_{k}^{x} )$ are both $(t-1)$-designs $(t\ge 2)$, then $(G,\mathcal{B}_{k}^{x} )$ is a $t$-design.
		\end{enumerate}
	\end{theorem}
	\begin{proof}
		When $(G,\mathcal{B}_{k}^{x} )$ is a $t$-design, then it is also a $1$-$(n,k,r)$ design with $r=b_{k}^{x}\cdot k/n<b_{k}^{x}$. Thus, $r_{k}^{x}(0)=r<b_{k}^{x}$, then $ \mathcal{B}_{k}^{x,*}\ne \emptyset$. For any $(t-1)$-subset $T=\{g_{1},\cdots, g_{t-1}\}$ of $G^{*}$, the blocks in $\mathcal{B}_{k}^{x} $ that contain the set $T$ consists of blocks that contain $0$ and blocks that do not contain $0$, where the blocks in the former family one-to-one corresponds to the blocks in $\mathcal{B}_{k-1}^{x,*}$ that contain $T$, and the blocks in the latter family  one-to-one corresponds to the blocks in $\mathcal{B}_{k-1}^{x,*}$ that contain $T$. Now, we have the following equation
		\begin{equation}\label{rt}
			r_{k}^{x}(T)=r_{k}^{x,*}(T)+r_{k-1}^{x,*}(T).
		\end{equation}
		If $(G,\mathcal{B}_{k}^{x} )$ is a $t$-design, then $(G,\mathcal{B}_{k}^{x} )$ is also a $(t-1)$-design, which follows that $r_{k}^{x}(T)$ is a constant independent of the choice of $T$. Moreover, for any $t$-subset $S=\{s_{1},\cdots,s_{t}\}\subseteq G$, the blocks in $\mathcal{B}_{k}^{x}$ that contain $S$ one-to-one correspond to the blocks that contain the $t$-set $\{0,s_{2}-s_{1},\cdots,s_{t}-s_{1}\}$ in $\mathcal{B}_{k}^{x}$, by the operation $-s_{1}$ in $G$ which induces a permutation in the family $\mathcal{B}_{k}^{x}$. This shows that $r_{k-1}^{x,*}(T)$ is also a constant independent of $T$, then by the Eq. (\ref{rt}), $r_{k}^{x,*}(T)$ is a constant, which proves $(i)$.\\ 
		If $(G^{*},\mathcal{B}_{k}^{x,*} )$ and $(G,\mathcal{B}_{k}^{x} )$ are both $(t-1)$-designs, $r_{k-1}^{x,*}(T)$ is a constant by Eq. (\ref{rt}), then by the one-to-one correspondence via the permutation mentioned earlier, $r_{k}^{x}(S)$ is a constant independent of the $t$-set $S$, then $(G,\mathcal{B}_{k}^{x} )$ is a $t$-design, which proves $(ii)$.
	\end{proof}
	\section{Applications to Elliptic Curve Codes Supporting $t$-Designs}\label{elliptic curve codes}
	The close relation between coding theory and $t$-design has attracted significant attention for many years. It is well known that a $t$-design may yield many linear codes and a linear code may support many $t$-designs, see \cite{D12020,D22020,D32020,T2020,Xiang2022,Xiang2023,Yan2024}. In this section, we characterize the support designs of minimum-weight codewords in some elliptic curve codes from $t$-designs held in their rational points groups. By such correspondence, we obtain a class of NMDS elliptic curve codes supporting $1$-designs, directly from the theories established in preceding section. Moreover, the weight distributions of the codes are given.\par               
	Let $\mathcal{C} $ be an $[n,k,d]$ linear code over the finite field $\mathbb{F}_{q}$. We index the coordinates of a codeword in $\mathcal{C} $ by $(1,2,\cdots, n)$ and define the set $\mathcal{P} (\mathcal{C} )=\{1,2,\cdots, n\}$. The \emph{support} of a codeword $\mathbf{c}=(c_{1},c_{2},\cdots, c_{n})$ is defined as $$\mathrm{Supp}(\mathbf{c})=\{i:c_i\neq0,i\in\mathcal{P}(\mathcal{C})\},$$ and by $\mathcal{H}_{d}(\mathcal{C} )$ we denote the set of distinct supports of all codewords with minimum Hamming weight $\mathrm{wt}(\mathbf{c})=d$.  We say that the minimum-weight codewords in the linear code $\mathcal{C}$ supports a $t$-design if the incidence structure $\left ( \mathcal{P} (\mathcal{C} ), \mathcal{H}_{d}(\mathcal{C} )\right ) $ is a $t$-$(n,d,\lambda)$ design. For further information about linear code and $t$-design, the reader is referred to \cite{D2022}. \par
	Now we introduce the Algebraic geometry (AG) codes and elliptic curve codes. AG codes are natural generalization of the Reed-Solomon codes, and the elliptic curve code is a type of AG codes, defined by algebraic curves of genus $g=1$, i.e, elliptic curves. We further assume that the characteristic of the finite field $\mathbb{F}_{q}$ is not $2$ or $3$, then an elliptic curves over $\mathbb{F}_{q}$ is given by the following equation
	$$y^{2}=x^{3}+ax+b,$$ where $a,b\in \mathbb{F}_{q} $. In addition, the discriminant of the curve $-16(4a^{3}+27b^{2})$ should be nonzero to ensure the smoothness of the curve. We denote the set of rational points on the elliptic curve $E$ by $E(\mathbb{F}_{q})$, which consists of the solutions of the equation and the infinity point $O$. The set $E(\mathbb{F}_{q})$ is a finite abelian group, with the zero element $O$. Moreover, the structure of the group $E(\mathbb{F}_{q})$ is either $\mathbb{Z}_n$ or $\mathbb{Z}_{n_1}\oplus\mathbb{Z}_{n_2}$ for some integers $n_{1},n_{2}$ with $n_{1}\mid n_{2}$. For more information about elliptic curves, see Silverman's book \cite{St2009}.\par
	Before introducing the definition of AG codes, we first fix the following notations
	\begin{itemize}
		\item $X$ is a smooth projective curve of genus $g$ over $\mathbb{F}_{q}$, and fix an algebraic closure $\overline{\mathbb{F}}_{q} $ of $\mathbb{F}_{q}$
		\item \(\overline{\mathbb{F}}_{q}(X)\) (resp. \(\mathbb{F}_{q}(X)\)) is the function field of \(X\) over \(\overline{\mathbb{F}}_{q}\) (resp. over \(\mathbb{F}_{q}\)).
		\item For any \(f \in \overline{\mathbb{F}}_{q}(X)^{\times}\), the principal divisor is \(\operatorname{div}(f) = \sum v_{Q}(f)Q\), where \(v_{Q}(f)\) is the valuation of \(f\) at \(Q\). 
		\item $X(\mathbb{F}_{q})$ is the set of $\mathbb{F}_{q}$-rational points on $X$.
		\item $\left \{ P_{1},P_{2},\cdots, P_{n} \right \} $ is proper subset of $X(\mathbb{F}_{q})$, and $D=P_{1}+P_{2}+\cdots+P_{n}$ is a divisor.
		\item $G$ is a divisor of degree $k$, where $2g-1\le k\le n-1$, and $\mathrm{Supp}(G)\cap \mathrm{Supp}(D)=\emptyset$.
		\item $\mathscr{L}(G)=\{f\in\mathbb{F}_{q}(X)^{\times} \mid (f)\ge -G \}\cup \{0\}$ is the Riemann-Roch space associated to $G$.
	\end{itemize}
	The reader is referred to \cite{S2009} for detailed information about AG codes.
	\begin{definition}
		The AG code $C_{\mathscr{L}}(D,G)$ is defined by the image of the evaluation map $\mathrm{ev}_{D}:\mathscr{L}(G)\longrightarrow \mathbb{F}_{q}^{n} $, given by 
		$$\mathrm{ev}_{D}:f\mapsto  \left ( f(P_{1}), f(P_{2}),\cdots,f(P_{n}) \right ) \in\mathbb{F}_{q}^{n}. $$
	\end{definition}
	By the well-known Riemann-Roch Theorem, the dimension $\dim C_{\mathscr{L}}(D,G)=k-g+1 $. The minimum distance of $C_{\mathscr{L}}(D,G)$ is lower bounded by $n-k$, for the number of zeros of functions in $\mathscr{L}(G)$ is upper bounded by $\deg G=k$, together with the Singleton bound, we have
	\begin{equation}\label{bound}
		n-k\le d\le n-k+g.
	\end{equation}
	When $X=E$, an elliptic curve over $\mathbb{F}_{q}$, Eq. (\ref{bound}) shows that the minimum distance of the elliptic curve code $C_{\mathscr{L}}(D,G)$ is either $n-k$ or $n-k+1$.The former case corresponds to an MDS code, while the latter defines an NMDS code, since the dual of an elliptic curve code is also an elliptic curve code \cite{S2009}. Very recently, Han and Ren \cite{HR2024} proved that the maximal length of $q$-ary MDS elliptic curve codes is close to $(1/2+\epsilon)\#E(\mathbb{F}_{q})$, which confirmed a conjecture proposed by Li, Wan and Zhang in \cite{LWZ2015}. In the following Proposition \ref{elli}, we introduce a general characterization of the MDS property from subset sums in the group $E(\mathbb{F}_{q})$.\par
	The following proposition connects $t$-designs from minimum weight codewords of elliptic curve codes and those from zero-sum subsets in the set $E(\mathbb{F}_{q})^{*}$. Note that part $(i)$ of the proposition is well-known \cite{HR2023,LWZ2015}. Denote by $\oplus$ the plus operator in the group $E(\mathbb{F}_{q})$.
	\begin{proposition}\label{elli}
		Let $E$ be an elliptic curve over $\mathbb{F}_{q}$, with $|E(\mathbb{F}_{q})|=n+1$. Let $D=P_{1}+P_{2}+\cdots+P_{n}$, where $\{ P_{1},P_{2},\cdots, P_{n}\} =E(\mathbb{F}_{q})^{*}$, and $G=kO$ $(1\le k\le n-1)$. Then we have 
		\begin{itemize}
			\item [(i)] The $[n,k,d]$ code $C_{\mathscr{L}}(D,G)$ is MDS, i.e., $d=n-k+1$, if and only if $\mathcal{B}_{k}^{*}=\emptyset$ in the group $E(\mathbb{F}_{q})$. Conversely, $C_{\mathscr{L}}(D,G)$ is NMDS, i.e, $d=n-k$, if and only if $\mathcal{B}_{k}^{*}\ne\emptyset$ in $E(\mathbb{F}_{q})$.
			\item [(ii)] If the code $C_{\mathscr{L}}(D,G)$ is NMDS, then the minimum-weight codewords support a $t$-design, i.e, $\left ( \mathcal{P} (\mathcal{C}_{\mathscr{L}}(D,G) ), \mathcal{H}_{n-k}(\mathcal{C}_{\mathscr{L}}(D,G) )\right ) $ is a $t$-design if and only if $\left ( E(\mathbb{F}_{q})^{*}, \mathcal{B}_{k}^{*}\right ) $ is a $t$-design.
			\item [(iii)] If $\left ( E(\mathbb{F}_{q}), \mathcal{B}_{k}\right )$ is a $1$-design, then the code $C_{\mathscr{L}}(D,G)$ is NMDS.
		\end{itemize}
	\end{proposition}
	\begin{proof}
		If $\mathcal{B}_{k}^{*}\ne\emptyset$, assume there are $\{P_{i_{1}},P_{i_{2}},\cdots, P_{i_{k}}\}\in E(\mathbb{F}_{q})^{*}$, such that in $E(\mathbb{F}_{q})$, \begin{equation}\label{sum}
			P_{i_{1}}\oplus P_{i_{2}}\oplus\cdots\oplus P_{i_{k}}=O.\end{equation}
		According to the isomorphism $E(\mathbb{F}_{q})\cong \mathrm{div}^{0}(E)/\mathrm{Prin}(\mathbb{F}_{q}(E))$, Eq. (\ref{sum}) is equivalent to saying that there is a unique function $f\in \mathbb{F}_{q}(E)$ (up to a constant), such that \begin{equation}\label{func}
			\mathrm{div}(f)=-kO+P_{i_{1}}+P_{i_{2}}+\cdots+P_{i_{k}}.
		\end{equation}
		By the fact that any principal divisor has a zero degree, Eq. (\ref{func}) is equivalent to saying that there exists some function $f\in \mathbb{F}_{q}(E) $, such that
		\begin{equation}
			\mathrm{div}(f)\ge -kO+P_{i_{1}}+P_{i_{2}}+\cdots+P_{i_{k}}.
		\end{equation}
		That is, there exists a function $f\in \mathscr{L}(G)$ with exactly $k$ zeros: $P_{i_{1}},\cdots, P_{i_{k}}$. This corresponds the existence of a codeword $c_{f}\in C_{\mathscr{L}}(D,G)$, where $c_{f}=ev_{D}(f)$, such that $\mathrm{Supp}(c_{f})=\left \{  1,2,\cdots, n \right \} \setminus\left \{ i_{1},i_{2},\cdots, i_{k}  \right \}$. Thus $d=n-k$, which proves $(i)$.\par
		To prove $(ii)$, equivalently, we prove $\left ( \mathcal{P} (\mathcal{C}_{\mathscr{L}}(D,G) ), \mathcal{H}_{n-k}(\mathcal{C}_{\mathscr{L}}(D,G) )\right ) $ is a $t$-design if and only if $\left ( E(\mathbb{F}_{q})^{*}, \mathcal{B}_{n-k}^{e,*}\right )$ is a $t$-design, where $e=\sum_{x\in E(\mathbb{F}_{q}) }^{} x$. To this end, consider the bijection $\Psi:\mathcal{B}_{n-k}^{e,*}\longrightarrow \mathcal{H}_{n-k}(\mathcal{C}_{\mathscr{L}}(D,G) ) $ given by 
		$$\Psi:\{( P_{j_{1}}, P_{j_{2}},\cdots,  P_{j_{n-k}})\mid P_{j_{1}}\oplus P_{j_{2}}\oplus\cdots\oplus P_{j_{n-k}}=e\}\mapsto  \{j_{1},j_{2},\cdots,j_{n-k}\}.$$  
		Indeed, given a support $\{j_{1},j_{2},\cdots,j_{n-k}\}$ of some codeword $c\in C_{\mathscr{L}}(D,G)$, as shown in the proof of $(i)$, there is a unique block in $\mathcal{B}_{k}^{*}$, with its complementary block  in $\mathcal{B}_{n-k}^{e,*}$. In addition, the bijection $\Psi$ naturally induces a preservation of incidence structures, this proves $(ii)$.\par 
		Finally, we have $\mathcal{B}_{k}^{*}\ne \emptyset$ if $\left ( E(\mathbb{F}_{q}), \mathcal{B}_{k}\right )$ is a $1$-$(n+1, k, r)$ design, where $r=b_{k}\cdot k/(n+1)<b_{k}$. Then the minimum distance $d=n-k$, which follows from $(i)$.
		
	\end{proof}
	In Proposition \ref{elli}, we build a bridge between the support designs in the elliptic curve code $\mathcal{C}_{\mathscr{L}}(D,G)$ (with $d=n-k$) and designs supported by in subset sums in $E(\mathbb{F}_{q})^{*}$. Further, when the rational points group $E(\mathbb{F}_{q})$ is the direct sum of two cyclic groups, that is, $E(\mathbb{F}_{q})\cong\mathbb{Z}_{n_1}\oplus\mathbb{Z}_{n_2}$, by Theorem \ref{t and t-1}, we immediately have the following characterization for such an elliptic curve code to support a $t$-design.
	\begin{theorem}\label{direct sum}
		Let $E$ be an elliptic curve over $\mathbb{F}_{q}$, with $E(\mathbb{F}_{q})\cong\mathbb{Z}_{n_1}\oplus\mathbb{Z}_{n_2}$, $1<n_{1}\mid n_{2}$, and let $n=n_{1}n_{2}-1$. Let the divisor $D=P_{1}+P_{2}+\cdots+P_{n}$, where $\{ P_{1},P_{2},\cdots, P_{n}\} =E(\mathbb{F}_{q})^{*}$, and $G=kO$ ($1\le k\le n-1$). Then in the $[n,k,n-k]$ code $\mathcal{C}_{\mathscr{L}}(D,G)$, the minimum-weight codewords support\label{key} a $t$-$(n,n-k,\lambda_{t})$ design, that is,  $\left ( \mathcal{P} (\mathcal{C}_{\mathscr{L}}(D,G) ), \mathcal{H}_{n-k}(\mathcal{C}_{\mathscr{L}}(D,G) )\right ) $ is a $t$-$(n,n-k,\lambda_{t})$ design if $n_{2}\mid k$ and $(\mathbb{Z}_{n_1}\oplus\mathbb{Z}_{n_2}, \mathcal{B}_{k})$ is a $(t+1)$-$(n+1, k, \lambda_{t+1} )$ design, where $t\ge 1$.
	\end{theorem}
	\begin{proof}
		By Proposition \ref{elli}, the incidence structure $\left ( \mathcal{P} (\mathcal{C}_{\mathscr{L}}(D,G) ), \mathcal{H}_{n-k}(\mathcal{C}_{\mathscr{L}}(D,G) )\right ) $ is a $t$-$(n,n-k,\lambda_{t})$ design if and only if $((\mathbb{Z}_{n_1}\oplus\mathbb{Z}_{n_2})^{*}, \mathcal{B}_{k}^{*})$ is a $t$-design, which holds if $(\mathbb{Z}_{n_1}\oplus\mathbb{Z}_{n_2}, \mathcal{B}_{k})$ is a $(t+1)$-$(n+1, k, \lambda_{t+1} )$ design, provided
		$n_{2}$ divides $k$, by Theorem \ref{t and t-1}.
	\end{proof}
	\begin{corollary}\label{direct sum of p}
		In Theorem \ref{direct sum}, if $n_{1}=n_{2}=p$, where $p$ is a prime, let $p$ divides $k$. Then in the $[p^{2}-1,k,p^{2}-k-1]$ code $\mathcal{C}_{\mathscr{L}}(D,G)$, the family of minimum-weight codewords supports a $1$-design.
	\end{corollary}
	\begin{proof}
		In this case, $\left ( \mathcal{P} (\mathcal{C}_{\mathscr{L}}(D,G) ), \mathcal{H}_{n-k}(\mathcal{C}_{\mathscr{L}}(D,G) )\right ) $ is a $1$-design if and only if $((\mathbb{Z}_{p}\oplus\mathbb{Z}_{p})^{*}, \mathcal{B}_{k}^{*})$ is a $1$-design. By Theorem \ref{t and t-1}, this is equivalent to saying that $(\mathbb{Z}_{p}\oplus\mathbb{Z}_{p}, \mathcal{B}_{k})$ is a $2$-design, which holds if and only if $p$ divides $k$, by Proposition \ref{x=0}.
	\end{proof}
	\begin{example}
		Let $E$ be the elliptic curve $y^{2}=x^{3}+3$ over $\mathbb{F}_{43}$, then the group $E(\mathbb{F}_{43})\cong \mathbb{Z}_{7}\oplus\mathbb{Z}_{7}$. Let the divisor $D=P_{1}+P_{2}+\cdots+P_{8}$, where $\{ P_{1},P_{2},\cdots, P_{8}\} =E(\mathbb{F}_{7})^{*}$, and let $G=kO$, where $7\mid k, k<49$. The minimum distance of the elliptic curve code $\mathcal{C}_{\mathscr{L}}(D,G)$ is $49-k$, because in the group $\mathbb{Z}_{7}\oplus\mathbb{Z}_{7}$, $b_{k}^{*}>0$ by Theorem \ref{thm 3c}. Then the $[49,k,49-k]$ ($7\mid k, k<49$) NMDS elliptic curve code $\mathcal{C}_{\mathscr{L}}(D,G)$ support a $1$-design, for $(\mathbb{Z}_{7}\oplus\mathbb{Z}_{7}, \mathcal{B}_{k}) $ is a $2$-design when $7\mid k$. This example has been verified by MAGMA.
	\end{example}
	
	According to the structure of rational point groups of elliptic curves, $E(\mathbb{F}_{q})$ is either $\mathbb{Z}_n$ or $\mathbb{Z}_{n_1}\oplus\mathbb{Z}_{n_2}$ for some integers $n_{1},n_{2}$ with $n_{1}\mid n_{2}$. Then the $t$-designs in the corresponding elliptic codes are linked to $t$-designs arising from the incidence structure  $(\mathbb{Z}_{n}^{*}  \mathcal{B}_{k}^{*})$ or  $((\mathbb{Z}_{n_{1}}\oplus\mathbb{Z}_{n_{2}})^{*}, \mathcal{B}_{k}^{*})$. As discussed in Section \ref{sec cyc}, preliminary results on $t$-designs in these incidence structures suggest a promising direction for further investigation.

	\section{Conclusion and Some Open Problems}\label{sec con} 
	This paper studies $t$-designs arising from subset sums in finite abelian groups, motivated by both their intrinsic mathematical interest and connections to coding theory. Our main contributions include:
	\begin{itemize}
		\item Characterization of necessary and sufficient conditions for $(G, \mathcal{B}_{k}^{x})$ to form a $1$-design, when $G$ is a finite abelian group of exponent $p^{m}$ and $pq$, and a conjecture regarding the non-existence of $2$-designs for non-elementary abelian $p$-groups.
		\item Some observations of $t$-design properties in cyclic groups and general non-cyclic abelian groups through the incidence structures $(G, \mathcal{B}_{k}^{x})$ and $(G^{*}, \mathcal{B}_{k}^{x,*})$, which applies to elliptic curve codes.
		\item Establishing connections between these combinatorial $t$-designs from subset sums and those from elliptic curve codes, demonstrating how our results yield NMDS elliptic curve codes supporting $1$-designs. We conclude by proposing two open problems for future research.
	\end{itemize}
	
	\begin{open problem}
		\rm{ If $G$ is the cyclic group $\mathbb{Z}_{n}$, under what conditions can the incidence structure $(\mathbb{Z}_{n},\mathcal{B}_{k}^{x})$ be a $1$-design? Similarly, can the incidence structure $(\mathbb{Z}_{n}^{*},\mathcal{B}_{k}^{x,*})$ be a $1$-design? If it can, try to determine the conditions.}
	\end{open problem}
	\begin{open problem}
		\rm{When $n_{1}$ and $n_{2}$ are both some power of a prime $p$, the conditions for $(\mathbb{Z}_{n_1}\oplus\mathbb{Z}_{n_2}, \mathcal{B}_{k}^{x})$ to be a $1$-design are derived in this paper. However, for general $n_{1}$ and $n_{2}$ with $n_{1}>1$ and $n_{1}\mid n_{2}$, the conditions are not yet known.}
	\end{open problem}
	The reader is invited to attack the two open problems and Conjecture \ref{conjecture}.

\end{document}